\documentclass[conference]{IEEEtran}
\IEEEoverridecommandlockouts
\usepackage{amsmath,amssymb,amsfonts}

\usepackage{algorithmic}
\usepackage{graphicx}
\usepackage{textcomp}
\usepackage{xcolor}
\usepackage[
    style=numeric-comp,
    sorting=none,
    backend=biber,
    sortcites=true,
    doi=false,
    firstinits=true,
    hyperref,
    isbn=false,
    eprint=false]{biblatex}

\renewbibmacro{in:}{}
\AtEveryBibitem{%
  \clearlist{language}%
  \clearfield{pages}%
}

\newcommand{\R}{\mathbb{R}}
\begin{document}
\title{Feedback Linearization of Car Dynamics for Racing via Reinforcement Learning\\
}

% \author{\IEEEauthorblockN{Chris Cai}
% \IEEEauthorblockA{\textit{Electrical Engineering and Computer Science} \\
% \textit{University of California, Berkeley}\\
% Berkeley, US \\
% xiangyu.cai@berkeley.edu}
% \and
% \IEEEauthorblockN{Michael Estrada}
% \IEEEauthorblockA{\textit{Electrical Engineering and Computer Science} \\
% \textit{University of California, Berkeley}\\
% Berkeley, US \\
% mkestrada2@berkeley.edu}
% \and
% \IEEEauthorblockN{Star Li}
% \IEEEauthorblockA{Electrical Engineering and Computer Science}\\
% \textit{University of California, Berkeley}\\
% Berkeley, US \\
% listar2000@berkeley.edu}

\author{\IEEEauthorblockN{Michael Estrada, Sida Li
, Xiangyu Cai}
\IEEEauthorblockA{\textit{Electrical Engineering and Computer Science}\\
\textit{University of California, Berkeley}\\
Berkeley, US\\
Email: \{mkestrada2, listar2000, xiangyu.cai\}@berkeley.edu}}

\maketitle

\begin{abstract}
Through the method of \textit{Learning Feedback Linearization}, we seek to learn a linearizing controller to simplify the process of controlling a car to race autonomously. A soft actor-critic approach is used to learn a decoupling matrix and drift vector that effectively correct for errors in a hand designed linearizng controller. The result is an exactly linearizing controller that can be used to enable the well developed theory of linear systems to design path planning and tracking schemes that are easy to implement and significantly less computationally demanding. To demonstrate the method of feedback linearization, it is first used to learn a simulated model whose exact structure is known, but varied from the initial controller, so as to introduce error. We further seek to apply this method to a system that introduces even more error in the form of a gym environment specifically designed for modelling the dynamics of car racing. To do so, we posit an extension to the method of learning feedback linearization; a neural network that is trained using supervised learning to convert the output of our linearizing controller to the required input for the racing environment. Our progress towards these goals is reported and next-steps in their accomplishment are discussed.
\end{abstract}


\section{Introduction}
The control of nonlinear systems is traditionally considered challenging because of the many, and varied classes of nonlinear systems that each display unique behaviors. typically this needs to be accounted for on a case-by-case basis. The field of geometric control seeks to address this by taking advantage of the structure of the system's dynamics to simplify related tasks like trajectory tracking and planning \cite{sastry2013nonlinear,isidori2013nonlinear}. The fundamental issue with this approach is the necessity for \textit{accurate} knowledge of the system's dynamics, which can involve a great deal of system identification work or otherwise simplifying the model of the dynamics. This is in contrast to the techniques of model free reinforcement learning,  \cite{bertsekas1996neuro, sutton2018reinforcement, sutton2000policy} which can be used to automatically compute optimal control strategies with no prior knowledge of the system's dynamics. Recent work in this direction has shown some promise but remains difficult in many cases due to poor sample complexity in many robotics applications \cite{schulman2015trust, lillicrap2015continuous, schulman2017proximal}. In this paper, we seek to demonstrate the technique of Learning Feedback Linearization (LFBL), first proposed in \cite{westenbroek2019feedback}. LFBL is designed to exploit the advantages of both geometric control and model-free reinforcement learning to optimize a geometric controller for a plant with unknown dynamics and limited structural assumptions about those dynamics. Specifically, it seeks to learn an optimal implementation of a technique called \textit{Feedback Linearization} (FBL). A successful feedback linearization scheme has the extremely useful property of rendering the input-output relationship of an otherwise \textit{nonlinear} plant to be \textit{linear}. This then enables the use of path planning and control techniques from the much more broadly applicable, and easily applied, field of linear systems theory \cite{ martin2003flat, kalman1960contributions, borrelli2017predictive}. While this technique has been used extensively in robotics, the fundamental concern with feedback linearization is the necessity for highly accurate models of the dynamics governing the robot \cite{grizzle2001asymptotically, ames2014rapidly, mellinger2011minimum, sastry2013nonlinear, isidori2013nonlinear,}. The field of adaptive control has seen extensive work on the construction of exactly linearizing controllers \cite{sastry2011adaptive, craig1987adaptive, sastry1989adaptive, nam1988model, kanellakopoulos1991systematic, umlauft2017feedback, chowdhary2014bayesian, chowdhary2013bayesian}, but this too has the limitation of needing highly structured representations of how and where the nonlinearities arise in the system's dynamics. The method of LFBL, in constrast, requires relatively limited knowledge about the structure of a nonlinear system to exactly linearize it.
\begin{figure}[!t]
    \centering
    \includegraphics[width=1.1\linewidth]{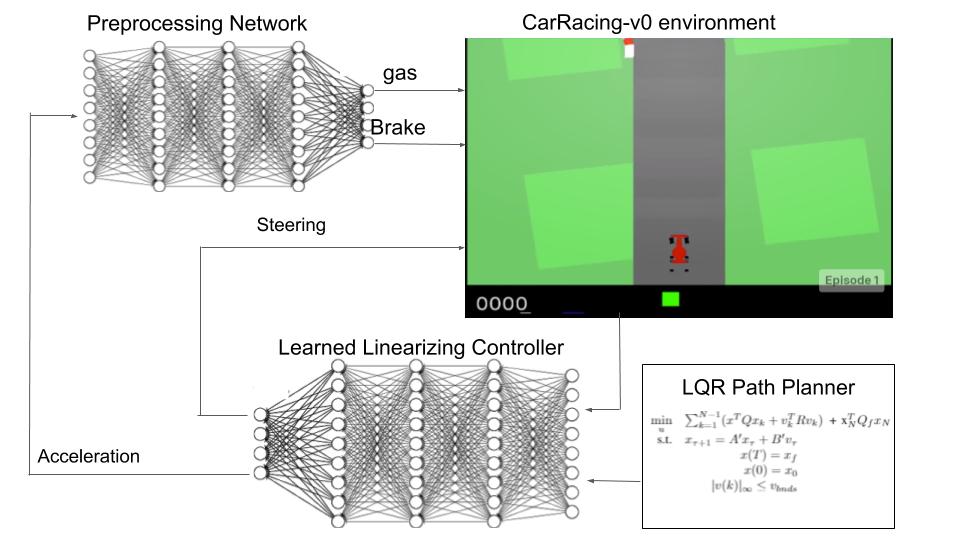}
    \caption{Diagram depicting the overall flow of the proposed control scheme.}
    \label{fig:flow chart}
\end{figure}
Beyond briefly presenting the method of FBL and extending it to LFBL, this paper is aimed at applying the method for the task of controlling 4-wheeled motor vehicles. Specifically, we want to employ this method to exploit the benefits of linear control in the development of algorithms for \textit{racing} cars. The nonlinear and nonholonomic nature of a car's dynamics, paired with their prevalence in our daily lives, makes the study of their efficient and effective control a tantalizing pursuit. A great deal of work has been done to accommodate the idiosyncrasies of path planning and control for vehicles, and autonomous racing has become a point of interest for profit, sport, and hobby alike \cite{DIYRoboCars,DeepRacer,F1th,RoboRace}. There is a mature body of literature on the design of Model Predictive Control (MPC) algorithms that seeks to create optimal paths based on an explicitly stated sets of dynamics and system constraints. There has also been some recent work towards learning various elements of MPC \cite{7487277,Borrelli_LMPC,Hewing}, With this background in mind, seek to use LFBL to learn the dynamics of a car. The resulting linear input-output relationship can be used to vastly simplify its planning and control. We demonstrate this by the construction of a simple path planner that uses a Linear Quadratic Regulator (LQR) structure paired with linear constraints to enforce the dynamics of the vehicle. In addition, we posit a novel extension to the method of LFBL; a preprocessing network that is designed to convert the values calculated by an exactly linearizing controller into an appropriate set of inputs for the system we wish to control. This conversion then allows for the application of a feedback linearizing controller to systems that do not share the same inputs as the linearizing controller.

The remainder of this paper is structured as follows: Section II: Theory, will briefly discuss the requisite theory for the method of FBL and how we formulate the problem of LFBL. Section III, Design Considerations, will address the specifics of the vehicle dynamics needed for this work as well as a variety of modifications and extensions made to the core methodology of LFBL to realize the desired controller. Section IV: Experementation, Results, and Discussion, will briefly address what experiments have been run to debug and validate the performance of our controller, and Section V: Conclusions and Future Work, will briefly summarize the goals and accomplishments of this work with an eye towards improvements that can be made moving forward.  
\section{Theory}

\subsection{Feedback Linearization} \label{FB}

To begin, we consider a system in the control affine form:
\begin{equation}
\label{eq:control_affine}
\begin{aligned}
    \dot{x} &= f(x) + g(x)u \\
    y &= h(x),
\end{aligned}
\end{equation}
The method presented here can address the "square" case of this form where the number of inputs matches the number of outputs. $x\in{\R^n}$ is the state, $u\in{\R^q}$ is the input, and $y\in{\R^q}$ is the output with the mappings $f:\R^n\rightarrow \R^n$, $g:\R^n\rightarrow \R^{nxq}$, and $h:\R^n\rightarrow \R^q$ are assumed smooth.

For the purposes of intuition, we will primarily address the single-input-single-output (SISO; q=1) case with the knowledge that the Multi-input Multi-output (MIMO) square case abides the same logic, albeit more mechanically obtuse. In the SISO case, we begin by taking time derivatives of $y=h(x)$ until the input $u$ appears and invert this relationship so-as to enforce the desired linear input-output behavior. Consider the first derivative of the output posed in (\ref{eq:control_affine}):

\begin{align*}\label{eq:Lie}
\dot{y} &= \frac{dh}{dx}(x) \cdot \Big(f(x) + g(x) u\Big)   \\
&= \frac{dh}{dx}(x) \cdot f(x) +\frac{dh}{dx}(x) \cdot g(x)u
\end{align*}

where $L_fh(x)=\frac{dh}{dx}(x)\cdot f(x)$ and $L_gh(x)=\frac{dh}{dx}(x)\cdot g(x)$ are called Lie derivatives and serve to quantify the rate of change of $y=h(x)$ along vector fields defined by $f$, and $g$ respectively. In the case where $L_gh(x) \neq 0$, $u$ appears in the derivative and we can exactly control the outputs $y=h(x)$ with a controller of the form:

\begin{equation}
\begin{array}{lcl} 
u(x,v) = \frac{1}{L_gh(x)}(-L_fh(x)+v) \\ 
\label{eq:controller}
\end{array}
\end{equation}

where $v \in R^q$ is the set of virtual inputs. This transformation between $u$ and $v$ serves to "cancel out" the nonlinearities of the original system and enforce the relationship $\dot{y}=v$. If, however, $L_gh(x)\equiv 0$, then the input can't directly affect the first derivative as desired and the control law (\ref{eq:controller}) is not well defined. In this case, we can continue taking derivatives of $y(x)$ until the input appears. When the input appears in the $\gamma^{th}$ derivative, the derivative can be expressed as:

\begin{equation}
\begin{array}{lcl} 
y^{(\gamma)}=L^{\gamma}_fh(x)+L_gL_f^{\gamma-1}h(x) u \\ 
\label{eq:Lie_higher_order}
\end{array}
\end{equation}
where $L^{\gamma}_fh(x)$ and $L_gL_f^{\gamma-1}h(x)$ are higher order Lie Derivatives. As before, if $L_gL_f^{\gamma-1}h(x)\neq 0$ then the control law can be rewritten as:

\begin{equation}
\begin{array}{lcl} 
u(x,v) = \frac{1}{L_gL_f^{\gamma-1}h(x)}(-L_f^{\gamma}h(x)+v) \\ 
\label{eq:controller_higher_order}
\end{array}
\end{equation}

As above, this choice of control law enforces the relationship $y^{\gamma}=v$, where gamma is called the "relative degree" of the nonlinear system.

This process can be extended readily to square, MIMO systems [$q > 1$]. By successively taking the derivatives of each output variable until at least one of the inputs appears, we can obtain an input-output relationship of the form:

\begin{equation}
\begin{array}{lcl} 
[y_1^{\gamma_1},...,y_p^{\gamma_q}]^T=b(x)+A(x)u \\ 
\label{eq:controller_higher_order}
\end{array}
\end{equation}

where the square matrix $A(x)\in\R^{q \times q}$ is called the "decoupling matrix" and $b(x) \in \R^q$ is called the "drift term". If $A(x)$ is invertable, then the control law can again be written as:
\begin{equation}
\begin{array}{lcl} 
u(x,v)= A^{-1}(x)(-b(x)+v)
\label{eq:MIMO_controller}
\end{array}
\end{equation}
with virtual inputs $v \in \R^{q}$ corresponding to the decoupled linear system in (\ref{eq:MIMO_controller}):
\begin{equation}
\begin{array}{lcl} 
[y_1^{\gamma_1},...,y_q^{\gamma_q}]^T=[v_1,....,v_q]^T
\label{eq:linear_system}
\end{array}
\end{equation}

where the $j^{th}$ entry of $v$ is denoted $v_j$, $y_j^{\gamma_j}$ is the $\gamma^{th}$ derivative of the $j^{th}$ output, and $\gamma = (\gamma_1,...,\gamma_q)$ is called the "relative degree" of the system. The system (\ref{eq:linear_system}) is Linear Time invariant (LTI) and can be represented in the "normal form":
with
\begin{equation}
\begin{array}{lcl} 
\dot{\xi_r} = A'\xi_r + B'v_r
\label{eq:normal_form}
\end{array}
\end{equation}
With state variables $\xi_r=(y_1,\dot{y_1},...,y_1^{\gamma_1-1},...,y_q,...,y_q^{\gamma_q-1})$ and $A' \in \R^{\gamma \times \gamma}$ and $B' \in \R^{\gamma \times q}$. This simplified normal form can be used to more easily create a desired trajectory for the output of the nonlinear system using the techniques of linear control theory.
A key insight from this analysis is the manner in which the transformed LTI system is made to behave. Specifically, the normal form is formulated in such a way that the $\gamma_j^{th}$ derivative of the $j^{th}$ output is the value being directly controlled while the actual output, $y_j$ is related to this value by a chain of integrators, e.g. we directly control the acceleration of a vehicle to achieve the desired position values.
\subsection{Learning Feedback Linearization}
The primary shortcoming of feedback linearization as a methodology is it's inherent dependence on \textit{exact} knowledge of the dynamics governing the plant, or, the system that is being controlled. This naturally leads to two sets of dynamics that must be considered in practice, the plant dynamics that represent the true behavior of the system, and the model dynamics designed by the engineer attempting to control the system; consider:
\begin{equation}
\label{eq:control_affine_plant}
\begin{aligned}
    \dot{x}_p &= f_p(x_p) + g_p(x)u_p \\
    y_p &= h_p(x),
\end{aligned}
\end{equation}

\begin{equation}
\label{eq:control_affine_model}
\begin{aligned}
    \dot{x}_m &= f_m(x_m) + g_m(x)u_m \\
    y_m &= h_m(x),
\end{aligned}
\end{equation}
Where a subscript $p$ denotes that we are referring to the plant and a subscript $m$ indicates that we are referring to the model. With the technique of Learning Feedback Linearization, we aim to rectify the aforementioned model mismatch by learning how to appropriately deviate the nonlinear input $u_p$ to better match the output of the plant to what we would expect based on the model. In this way, we are not only learning an exactly linearizing controller, but also forcing the plant to behave according to the dyanamics designated by our model.
For learning feedback linearization to be applicable, we first require that the system(s) (\ref{eq:control_affine_plant}), (\ref{eq:control_affine_model}) have the same, well-defined relative degree:
\begin{equation}
\label{eq:pm_rel_degree}
\begin{aligned}
    (\gamma_1^p,...,\gamma_q^p) = (\gamma_1^m,...,\gamma_q^m)
\end{aligned}
\end{equation}
With this assumption, we can guarantee the existence of feedback control laws of the form:
\begin{equation}
\begin{array}{lcl} 
u_p(x,v)= \beta_p(x) + \alpha_p(x)v
\label{eq:model_controller}
\end{array}
\end{equation}
\begin{equation}
\begin{array}{lcl} 
u_m(x,v)= \beta_m(x) + \alpha_m(x)v
\label{eq:plant_controller}
\end{array}
\end{equation}
With $\beta_p(x), \beta_m(x) \in \R^{q}$, and $\alpha_p(x),\alpha_m(x) \in \R^{q \times q}$. $u_p$ is, as of yet, unknown, but owing to the structure and assumptions developed here, we can relate the dynamics of the plant and model as follows:
\begin{equation}
\begin{array}{lcl} 
\beta_p(x) = \beta_m(x) + \Delta\beta(x)
\label{eq:beta_controller}
\end{array}
\end{equation}
\begin{equation}
\begin{array}{lcl} 
\alpha_p(x) = \alpha_m(x) + \Delta\alpha(x)
\label{eq:alpha_controller}
\end{array}
\end{equation}

Here, $\Delta\beta(x)$ and $\Delta\alpha(x)$ are two continuous functions with parameterized estimates:

\begin{equation}
\begin{array}{lcl} 
\Delta\beta(x) \approx \beta_{
\theta_1}(x), \Delta\alpha(x) \approx \alpha_{\theta_2}(x)    
\label{eq:approx_controller}
\end{array}
\end{equation}

Where $\theta_1, \theta_2$ represent the parameters to be learned through experimentation of the reinforcement learning agent. It is through learning the approximations (\ref{eq:approx_controller}) that the model mismatch is rectified and exact feedback linearization of the dynamic system in question can be achieved with the control law:
\begin{equation}
\begin{array}{lcl} 
\hat{u}_{\theta}(x,v)=[\beta_m(x)+\beta_{\theta_1}] + [\alpha_m(x) + \alpha_{\theta_2}]v  
\label{eq:uhat_controller}
\end{array}
\end{equation}

\subsection{Formulating the Learning Problem}
\subsubsection{continuous time}
We know from section \ref{FB} that the dynamics relating the inputs to the outputs under control law (\ref{eq:uhat_controller}) are of the form: 
\begin{align}{} 
y^{(\gamma)}=\underbrace{b_p(x) + A_p(x)\hat{u}_{\theta}(x.v)}_{W_\theta(x,v)}
\label{eq:dynamics_relation}
\end{align}
Where $b_p(x)$ and $A_p(x)$ are unknowns. We seek to find the parameters $\theta$ such that $W_{\theta^*}(x,v)\approx v$, leading to the choice of a point-wise loss $\ell:\R^n\times \R^q \times \R^{K_1+K_2} \rightarrow \R$:
\begin{align}{} 
\ell(x,v,\theta)=||v-W_{\theta^*}(x,v)||_2^2
\label{eq:continuous_loss}
\end{align}
This relatively simple 2-norm loss measures the learned controller's performance in linearizing the plant dynamics as a function of the control $\hat{u}_{\theta}$ and the state $x$ when virtual input $v$ is applied to the linear reference model whose behavior we seek to match. Next we can define a probablity distribution $X$ over $\R^n$ with $V$, the uniform distribution over the set $\{v\in \R^q:||v||\leq1\}$ and define the weighted loss:
\begin{equation} 
L(\theta)=\mathbb{E}_{x \sim X,v \sim V}[ l(x,v,\theta)]
\label{eq:weighted_continuous_loss}
\end{equation}
With this, we choose to formulate our learned controller to optimize over parameters, $\theta$ to solve the problem:
\begin{equation}
\min_{\theta}L(\theta)
\label{eq:weighted_continuous_opt}
\end{equation}

The primary concern in formulating the learning problem like this is the unknown nature of the terms that constitute $W_{\theta^*}(x,v)$, however, this can be circumvented by measuring the outputs of the plant, $y^{(\gamma)}$ as there is an equivalence between the two as defined in (\ref{eq:dynamics_relation}).
Further detail on the continuous time formulation of this problem, including the conditions under which global minima can be found, can be found in \cite{westenbroek2019feedback}.
\subsubsection{discrete time}
While the theory above is developed generally in the continuous-time case, it is of practical use for many systems, including the work presented here, to formulate the learning problem in discrete-time.
 To do so, we first define $\Delta t>0$ to denote the sampling rate of our system and $t_k=k\Delta t$ for each $k \in \mathbb{N}$. Then, defining $x(\cdot)$  as the state trajectory of the plant with $x_k=x(t_k)$ and let $u_k$ denote the control applied on interval $[t_k,t_{k+1})$. With these, we can setup the difference equation for the plant dynamics:
 \begin{equation}
x_{k+1} = x_k +\underbrace{\int_{t_k}^{t_{k+1}}{f_p(x(t))+g_p(x(t))u_kd\tau}}_{F(x_k,u_k)}
\label{eq:disc_x_dyn}
\end{equation}
Now, letting $\xi(\cdot) = (y_1(\cdot),...,y_1^{(\gamma_1-1)}(\cdot),...,y_q(\cdot),...,y_q^{(\gamma_q-1)}(\cdot))$ represent the linearized system state and defining $\xi_k = \xi(t_k)$, we can integrate (\ref{eq:dynamics_relation}) over $[t_k, t_{k+1})$ and obtain a new difference equation in terms of the linearized system states:
 \begin{equation}
\xi_{k+1} = e^{A \Delta t}\xi_k + H(x_k,u_k)
\label{eq:disc_xi_dyn}
\end{equation}
here, $H(x_k,u_k)$ is, in general, no longer affine in $u_k$, presenting the main issue in applying the previously developed continuous-time theory here. Nonetheless, we can still try to enforce an approximately linear relationship between a virtual input $v_k$ and the subsequent outputs. for ease of Notation, define $\bar{A}=e^{A\Delta t}$, and $\bar{B}=\int^{\Delta t}_0e^{At}Bdt$, then the normal form (\ref{eq:normal_form}) can be rewritten as:
 \begin{equation}
\xi_{k+1} = \bar{A}\xi_k + \bar{B}v_k
\label{eq:disc_xi_dyn2}
\end{equation}
This equation defines the linear relationship we would like to learn and enforce by our choice of control $u_k=\hat{u}_{\theta}(x_k,v_k)$. With this being the case, we define our discrete time point-wise loss:
\begin{equation}
\bar{\ell}(x_k,u_k,v_k)= ||\bar{B}v_k - H(x_k,u_k)||_2^2
\label{eq:weighted_disc_loss}
\end{equation}
 As above in the continuous case, this value can be readily calculated through observation of the system's outputs and numerically differentiating their values between time-steps, this loss metric is designed to measure how well the control $u_k$ truly corrects for the unknown components and nonlinearities of the plant.
We can then define the discrete time reinforcement learning problem over the parameters of $\hat{u}_{\theta}$ as:
\begin{equation}
\begin{array}{rrclcl}
\displaystyle \min_{\theta} & \multicolumn{3}{l}{\mathbb{E}_{x_0 \sim X,v_k \sim V,w_k \sim \mathcal{N}(0,\sigma_w^2)}[\sum_{k=1}^N\bar{\ell}(x_k,u_k,v_k)]}\\
\textrm{s.t.} &  x_{k+1} = x_k + F(x_k,u_k)\\
& x_0=x_0   \\
& u_k = \hat{u}_{\theta}(x_k,v_k) +w_k
\end{array}
\label{eq:DLQR}
\end{equation}
Where $N$ is the length of the training episodes, and $k$ represents the current time step. In this formulation, initial conditions are sampled using the desired state distribution $X$, and inputs are sampled from $V$ at each time step. the additional noise term, $w_k$ serves to encourage exploration. With this formulation, local optima can be found using a variety of reinforcement learning algorithms \cite{SuttonBarto,TRPO,PPO,DDPO} in tandem with with hardware experiments to evaluate the loss function. 
\section {Design Considerations}
\subsection{Feedback Linearization of Vehicle Dynamics}
For the purposes of this work, we address only the lateral dynamics of a vehicle, or, the dynamics from a "top down" view. The intricacies of a vehicle's suspension, the possibility of tire slip, and the dynamics of the car's drive train will ideally be learned by the preprocessing network posited in section \ref{ppnet}. With this in mind, we consider a dynamically extended bicycle model of the vehicle dynamics:
\begin{equation}
\begin{array}{lcl} 
\dot{x}=Vcos(\psi + \beta) \\ 
\dot{y}=Vsin(\psi + \beta) \\
\dot{\psi}=\frac{V}{l_r}sin(\beta) \\
\dot{V}=a \\
\dot{\beta}=b
\label{eq:Bike}
\end{array}
\end{equation}

\begin{figure}[!b]
    \centering
    \includegraphics[width=1.1\linewidth]{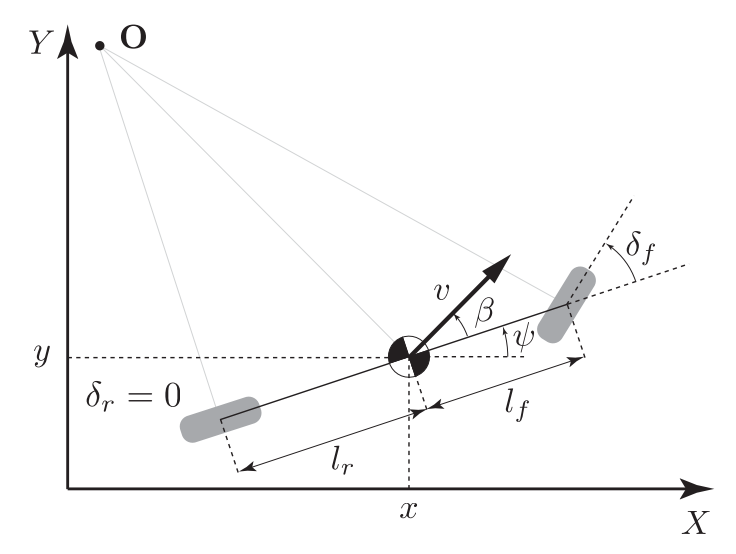}
    \caption{Visual representation of kinematic bicycle model}
    \label{fig:bike}
\end{figure}

Where x and y are the coordinates of the vehicle's center of mass in an inertial frame $(X,Y)$, $\psi$ is the inertial heading of the car, $V$ is the velocity of the car, $l_r$ is the distance from the center of mass to the vehicle's rear axle, and $\beta$ is the angle of the velocity of the center of mass with respect to the longitudinal axis of the car. Here, we choose $x$ and $y$ as the outputs of the system that we seek to control and further choose $a$ and $b$ as the inputs through which we will do so. A visualization of this system is provided below in figure \ref{fig:bike}. A first-principles analysis of this system reveals further:

\begin{equation}
\beta = tan^{-1}(\frac{l_r}{l_f+l_r}(tan(\delta_f))
\label{eq:Beta}
\end{equation}

where $\delta_f$ is the angle of steering of the front tires and $l_f$ is the distance from the center of mass to the vehicle's front axle. While $\delta_f$ is more intuitively meaningful as a choice of input, the complicated nature of its relationship to the vehicle's dynamics make it preferable for our formulation to control $\dot{\beta}$ directly and back-calculate $\delta_f$ uniquely in the case of a steering angle limited to at most $\pm 90$ degrees from a neutral orientation. This easily covers the range of any standard vehicle's max steering angle. 

With this form of the dynamics, we can start deriving the system's linearized normal form. We first take derivatives of the outputs with respect to time until the inputs appear:

\begin{equation}
\begin{array}{lcl} 
x = x \\
\dot{x} = Vcos(\psi + \beta) \\ 
\ddot{x} = \dot{V}cos(\psi + \beta) - \dot{\beta}Vsin(\psi + \beta)-\frac{V}{l_r}sin(\psi + \beta) \\
\ddot{x} =  cos(\psi + \beta)a - Vsin(\psi + \beta)b - \frac{V}{l_r}sin(\psi + \beta)
\label{eq:x_derivs}
\end{array}
\end{equation}

similarly, the derivatives of y can be found:

\begin{equation}
\begin{array}{lcl} 
y = y \\
\dot{y} = Vsin(\psi + \beta) \\ 
\ddot{y} = \dot{V}sin(\psi + \beta) + \dot{\beta}Vcos(\psi + \beta) +\frac{V}{l_r}sin(\psi + \beta) \\
\ddot{y} = sin(\psi + \beta)a + Vcos(\psi + \beta)b +\frac{V}{l_r}sin(\psi + \beta)
\label{eq:y_derivs}
\end{array}
\end{equation}

With this, our vector relative degree is:

\begin{equation}
\left[
\begin{array}{ccc}
\gamma_1 \\
\gamma_2 
\end{array}
\right] =
\left[
\begin{array}{ccc}
2 \\
2
\end{array}
\right]
\label{eq:rel_deg}
\end{equation}
By applying the transformation $x=\xi_1^1$, and $y=\xi_2^1$, the normal form can be written as: 

\begin{equation}
\begin{array}{lcl} 
\dot{\xi}_1^1 = \xi_2^1  \\
\dot{\xi}_2^1 = \xi_3^1  \\
\dot{\xi}_3^1 =  cos(\psi + \beta)a - Vsin(\psi + \beta)b -\frac{V}{l_r}sin(\psi + \beta)\\ 
\dot{\xi}_1^2 = \xi_2^2  \\
\dot{\xi}_2^2 = \xi_3^2  \\
\dot{\xi}_3^2 =  sin(\psi + \beta)a + Vcos(\psi + \beta)b+\frac{V}{l_r}sin(\psi + \beta) \\ 

\label{eq:normal_form_car}
\end{array}
\end{equation}
with corresponding decoupling matrix:
\begin{equation}
A=
\left[
\begin{array}{ccc}
 cos(\psi + \beta) & -Vsin(\psi + \beta)\\
 sin(\psi + \beta) & Vcos(\psi + \beta)
\end{array}
\right]
\label{eq:decoupling}
\end{equation}
and associated drift term:
\begin{equation}
b=
\left[
\begin{array}{ccc}
-\frac{V}{l_r}sin(\psi + \beta) \\
\frac{V}{l_r}sin(\psi + \beta) 
\end{array}
\right]
\label{eq:decoupling}
\end{equation}

This yields the equations for our linearizing inputs:

\begin{equation}
\begin{array}{ccc}
u= -A^{-1}b + -A^{-1}v \\
\end{array}
\label{eq:inputs}
\end{equation}

where $v$ is the virtual input to the desired linear system and $u$ is the corrected input that accounts for the nonlinearities of the system. Finally, we arrive at the linearized dynamics relating the virtual inputs to the outputs:
\begin{equation}
\begin{array}{ccc}
\ddot{x} = v_1  \\
\ddot{y} = v_2
\end{array}
\label{eq:linearized system}
\end{equation}

It is clear from inspection that the appropriate choice of $v_1$ and $v_2$ allows us to directly control the desired output, $x$ and $y$. The resulting \textit{linear} system as expressed in (\ref{eq:normal_form}) is

\begin{equation}
A'=
\left[
\begin{array}{cccccc}
0 & 1 & 0 & 0  \\
0 & 0 & 0 & 0 \\
0 & 0 & 0 & 1 \\
0 & 0 & 0 & 0 \\
\end{array}
\right], 
B'=
\left[
\begin{array}{ccc}
0 & 0 \\
1 & 0 \\
0 & 0 \\
0 & 1 \\
\end{array}
\right]
\label{eq:A'}
\end{equation}
\begin{equation}
C'=
\left[
\begin{array}{cccccc}
1 & 0 & 0 & 0 \\
0 & 0 & 1 & 0 \\
\end{array}
\right]
\label{eq:C'}
\end{equation}
\subsection{Path Planning and Tracking}\label{planning}
The transformation of the previously nonlinear vehicle dynamics into a linear system of the form (\ref{eq:normal_form}) with dynamics governed by (\ref{eq:A'}), (\ref{eq:C'}) enables the use of a wide range of possible linear control and path planning techniques that present greater opportunity for guarantees in performance and safety. Central to the path tracking approach used in this work is the continuous-time LQR: 

\begin{equation}
\begin{array}{rrclcl}
\displaystyle \min_{u} & \multicolumn{3}{l}{\int_{0}^{\infty}(x^{T}Qx+v^{T}Rv)}\\
\end{array}
\label{eq:LQR}
\end{equation}
Where $x\in \R^{\gamma}$ is the system state, $v \in \R^{q}$ is the input and $Q\in\R^{\gamma \times \gamma}$ and $R\in\R^{q \times q}$ are weight matrices determined by designer. The solution of this optimization problem can be readily found by solving the continuous-time algebraic ricatti equation:
\begin{equation}
\begin{array}{rrclcl}
A^TP+PA-PBR^{-1}B^TP+Q=0
\end{array}
\label{eq:CARE}
\end{equation}
Solving this equation for a unique, positive definite solution, P, provides sufficient information to construct an optimal feedback control law:
\begin{equation}
\begin{array}{rrclcl}
u^*=Fx = R^{-1}B^TPx
\end{array}
\label{eq:opt_law}
\end{equation}

The ease of solution of equation (\ref{eq:CARE}), paired with the simple control law (\ref{eq:opt_law}) and guarantees of optimality in control make this ideal as a feedback control strategy to regulate the car to the commands of a higher-level path-planner, which is itself a discrete form of the LQR Controller: 
\begin{equation}
\begin{aligned}
\min_{u} & \sum_{k=1}^{N-1}(x_k^{T}Qx_{k}+v_{k}^{T}Rv_{k}) + x_{N}^TQ_fx_N\\
\textrm{s.t.} \quad & x_{\tau+1}=A'x_{\tau}+B'v_{\tau}\\
& x(T) = x_f   \\
& x(0) = x_0   \\
& |v(k)|_{\infty} \leq v_{bnds}
\end{aligned}
\label{eq:DLQR}
\end{equation}
Here, $N$ is the number of execution steps in the desired look-ahead time, $T$. $\tau$ represents the current time-step, $x_f$ and $x_0$ represent the desired final and initial locations of the vehicles \textit{linear} state respectively, $x_{\tau}$ and $v_{\tau}$ represent the linear system state and input, respectively, at time-step $\tau$. Finally, $Q \in \R^{\gamma \times \gamma}$, $Q_f\in \R^{\gamma \times \gamma}$, and $R\in \R^{q \times q}$ are weight matrices determined by the designer. Further,
\begin{equation}
\begin{aligned}
v=[v_1,..,v_{N-1}] \\
x=[x_1,...,x_N]
\end{aligned}
\label{eq:x and u}
\end{equation}
This configuration is ideal for multiple reasons. While the unconstrained version of (\ref{eq:DLQR}) can be solved readily using a ricatti equation similar to (\ref{eq:CARE}), we choose to formulate this system with the additional constraints of a starting and ending location, Other constraints are added to enforce mechanical limitations like only having the ability to accelerate so quickly. This effectively forces the planner to account for the boundary conditions of the car's starting location, mechanical limitations, and where it should be at the end of its look-ahead period. The addition of constraints to enforce waypoint tracking between the start and finish locations can be used to further constrain the planner to a desired path. In the exactly linearizing case, both the state path and inputs that we calculate in this approach can be used limiting the need for the continuous time LQR solution (\ref{eq:opt_law}), However, in the approximate controller case, we will need a more active form of feedback control to correct for the unknown true dynamics of the system. It is convenient to note that in both cases, (\ref{eq:DLQR}) a linearly constrained quadratic programming problem that is easily solvable by most off-the-shelf solvers quickly and thus is ideal for a waypoint following-based path planning approach. Specifically, this formulation can be iteratively solved as the system evolves to update the planned path as the vehicle successively approaches and passes new waypoints, akin to an MPC approach.

\subsection{Preprocessing Network}\label{ppnet}
The gym environment `CarRacing-v0` environment takes in tuples of actions in the form $\{s, g, b\}$, corresponding to the steering, gas, and brake applied respectively. In an end-to-end reinforcement learning task, there exists no need to investigate the detailed effects and properties of individual actions, as the learning algorithm would implicitly learn some representations of the actions, such as a Q functions, through the training process \cite{DQN}. 

In this task, however, the gym environment only offers the vehicle dynamics that the LFBL controller is applied to. Since the gym environments original reward function is unused, the learning algorithm cannot capture the action representations as stated above. Additionally, according to the bicycle model, the LFBL seeks to control the steering angle and the acceleration directly on the vehicle, but the car racing dynamics require gas and brake parameters instead. The mismatch in control inputs makes it necessary to train a preprocessing network with supervised learning.

Since relationship between linear accelearation and the control inputs is not too complicated (see Figure 1), a shallow feed-forward neural network can approximate the function well. The input layer has size $(1 \times 200)$ followed by a leakyReLU activation, and the output has size $(200 \times 1)$ with identity output activation. The Adam optimizer is used with a learning rate of 0.05 over the MSE loss. 

\subsection{Model-free Deep RL Training}
We utilize deep reinforcement learning to train our LFBL controller since neural networks are capable of learning and representing the complex relationship in nonlinear control systems. Since the accurate true dynamics of our system is unknown, we further restrict our attention to Soft Actor-Critic (SAC), a state-of-the-art model-free deep RL algorithm. SAC is an off-policy algorithm that optimizes over a stochastic policy with entropy regularization, which augments its objective to be
\begin{equation}
    J(\pi) = \sum_{t = 0}^T \mathbb{E}_{(s_t, a_t) \sim \rho_\pi}[r(s_t, a_t) + \alpha H(\pi(\cdot|s_t))]
\end{equation}
where $H(\pi(\cdot|s_t))$ is the entropy term that measures the randomness of the policy under $s_t$ and $\alpha$ is a temperature hyper-parameter that determines the relative importance of the entropy \cite{SAC}. Detailed implementations and novelty designs \cite{SACsummary} of SAC will not be mentioned in this work. The main reasons that we use SAC in training our linearizing controller are (1) it attains good balance between performance and sample efficiencey and (2) it is in-built in the OpenAI Spinning Up library. 

\subsection{Overall Control Architecture}
To briefly summarize the overall control architecture posited in this paper, we direct the reader to figure \ref{fig:flow chart} which is a visual representation of how the path and corresponding control signals are created and fed to the environment. In general, the entire system is guided by a path planner that calculates an optimal path from a desired starting point to a desired end-point, as-in (\ref{eq:DLQR}). Then the LFBL controller can be used to calculate an input to the nonlinear dynamics of the car that will behave as-in the desired linear relationship, (\ref{eq:linear_system}). The output of this controller, acceleration, is then converted into gas and brake signals that are fed directly to the racing environment in tandem with the other output of the linearizing controller, steering angle. The effect of these inputs is observed and used to further train the LFBL controller. Note that, in the case of the exact learned controller desribed below, the preprocessing network is not needed, as the dynamics being learned is strictly of the form of the LFBL controller. 

\section{Experimentation, Results, and Discussion}
\subsection{Exact Learned Controller}
We first sought to learn a controller that starts with a linearizing controller based on a set of nominal dynamics that can be exactly linearized relative to the true dynamics being learned. The setup for running experiments to this end involves the instantiation of two sets of dynamic equations of the form (\ref{eq:Bike}). The key difference is the differing choice of $l_r$. This inherently affects how the inertial heading evolves and has downstream consequences for how the outputs we seek to control evolve. In particular, this mostly effects the drift term. Thus the problem is for the agent to learn a \textit{corrective} drift term that accurately accounts for the difference in $l_r$ and causes the true dynamics to behave \textit{exactly} as the nominal dynamics would dictate. For the tests performed, the nominal dynamics were instantiated with an $l_r$ value of .5, idicating a distance between the front wheels and the center of mass of .5 meters.

With this setup we began training. The LQR planner discussed in section \ref{planning} was configured to plan a path for the car from a starting location of $[x,y]=[0,0]$ at $t=0$ to a final location of $[x,y]=[5,5]$ at $t=5$s with a time step of .02s, resulting in episodes of 249 steps of the environment. because the car begins aligned with the y axis, this involves turning, which will in turn require that the agent learn to correct for turning as well as acceleration. The step reward is calculated as the 2-norm distance of the vehicle's \textit{linear} state from the desired state calculated by the planner at that time. the sum of these values over an episode is therefore a good indication of the ability of the controller to follow the desired path and is consequently a good proxy for the quality of the learned controller.
\begin{figure}[!b]
    \centering
    \includegraphics[width=1.1\linewidth]{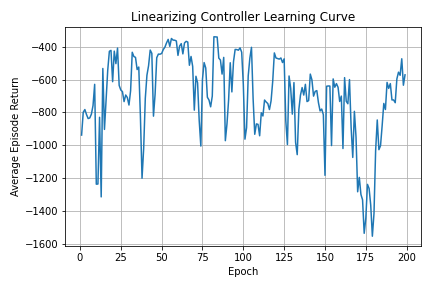}
    \caption{Plot comparing the path we want the agent to follow, the agent's true path and the path of the nominal linearizing controller with no learning}
    \label{fig:training}
\end{figure}
To begin, we removed the ability of the agent to actively effect the controller in order to set a baseline that reflects the performance of the controller as it was nominally designed with $lr=.5$. We found that in an episode of 249 steps the total loss for an episode was constant, as there was no way for the agent to interact with, and correct, the controller. The baseline average episode return value we observed was $-1.4 \times 10^{3}$. We next tried to run the algorithm using the default parameters of the SAC implementation being used and received a max average episode return of -350 at epoch 56 of 200. This is slightly less than a four-fold reduction, indicating significant improvements in the performance of the learned controller, a plot for the learning curve of this agent, as well as a visualization of the path it followed vs the path it was intended to follow are in figures \ref{fig:training} and \ref{fig:performance}, respectively. Once the performance of a naively tuned agent was found, we began a hyperparameter sweep to determine the best performing set of parameters. We opted to run these sweeps over a single variable, keeping the rest of the hyperparameters as the default setting used in the previous experiment. Together, we ran tests varying learning rate, the batch size of training updates, the discount factor, $\gamma$, as well as the polyak iterpolation factor, $\rho$, the result of most of these experiments was the same; within roughly 10 eopchs of starting the experiment, the agent would converge to an Average episode return fluctuating between -700 and -900, although it would occassionally show chaotic behavior where it either performed exceptionally well, as in the trial reported above, or blowing up and returning NaN values after several epochs. This seems to indicate that the performance of the agent is sensitive to the seed used in the training process, which was selected at random for each trial. We further note that these trials were run in roughly one hour per trial on both a personal laptop and desktop.

From inspection of \ref{fig:performance}, a few things are clear. While the desired path is relatively smooth, the path that is followed by the nominal controller with no learning performs quite poorly, diverging from the initial path very quickly, performing chaotic maneuvers and not arriving at the correct final location. This is to be expected as the controller is designed not to match the model that is actually being controlled. More interesting is the path followed by the learned controller, reported above to have an average episode return of -350. It is clear that this controller does not follow the path very closely, makes erroneous maneuvers,
and fails to reach the specified endpoint. It does, however, reach roughly the correct y-coordinate, $y=5$ and manages to stay much closer to the desired path than the the controller with no learning. There might be a few reasons for this lackluster performance. As previously mentioned, this trial appears to be the product of a good initial seed without any special hyperparameter tuning. Thus, performing a second hyperparameter sweep using this initial seed will likely improve performance significantly.  
\begin{figure}[!b]
    \centering
    \includegraphics[width=1.0\linewidth]{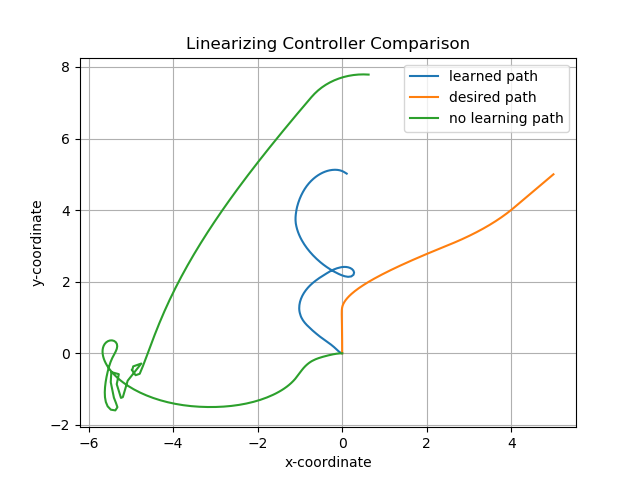}
    \caption{Plot comparing the path we want the agent to follow, the agent's true path and the path of the nominal linearizing controller with no learning}
    \label{fig:performance}
\end{figure}
\subsection{Approximate Learned Controller}
The next case to consider is the approximate learned controller. here the true dynamics are modeled by the car racing environment discussed in section III. Beyond simply learning a corrective term for \ref{eq:Bike}, this set of true dynamics also models more complicated behaviors relating to the application of brakes, steering, and gas as well as accounting for the delays in their effect on the acceleration and steering that we seek to control. This is where the preprocessing network is used to learn those more complicated relationships and relate them to the acceleration and steering angle that we command from our linearized controller. This involves not only a well-trained linearizing controller, but a well trained preprocessing network. As such, we though it best to decouple the problem of training the controller and preprocessing networks before we could begin experiments to this end in-earnest. While an initial version of the preprocessing network has been trained as per section \ref{ppnet_train}, the fact that we are still in the process of training an optimally performant linearizing controller has limited our ability to test the performance of the controller and pre-processing network combined.

\subsection{Preprocessing Network Training}\label{ppnet_train}
The aim of the preprocessing network is to take in the desirable linear acceleration that the feedback linearizer outputs, and return the tuple $\{gas, brake\}$ which can be directly sent as a control action to the gym environment. 

\begin{figure}[!b]
    \centering
    \includegraphics[width=1.1\linewidth]{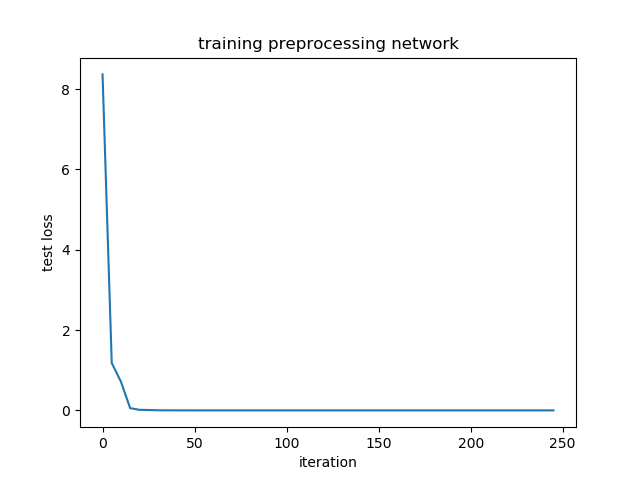}
    \caption{Plot depicting the learning curve of the preprocessing network.}
    \label{fig:preprocess}
\end{figure}

The training data set is collected by sampling random actions from the `CarRacing-v0` environment with the steering angle fixed to be zero. Once an action vector is sampled, it is fixed and kept applied to the environment for certain number of steps (usually 10). The acceleration is calculated by taking the average of the differences in velocity between consecutive steps. Our final data set consists of 5000 $\{acceleration, gas, break\}$ observations, where the last two variables serve as the label in the supervised learning task. A plot of the learning curve for the initial preprocessing network implementation is shown in figure \ref{fig:preprocess}. While it appears to have convereged by iteration 20 of training, we have not yet had the opportunity to  test it in tandem with the linearizing controller.
\section{Conclusions and Future Work}
While traditional methods of geometric control can drastically simplify the control of complex nonlinear systems, they often require a great deal of knowledge about the structure and parameters that govern the dynamics of the system to be controlled. In stark contrast to this, model free methods of reinforcement learning offer the ability to learn controllers for dynamic systems with almost no prior knowledge of the system's dynamics, but at the cost of requiring large data sets and a lot of computer power. The method of learning feedback linearization is devised by \cite{westenbroek2019feedback} to fuse these two sets of techniques to effectively benefit from the advantages of both. A feedback linearizing controller is designed using a nominal set of dynamics, and model mismatch between the plant and controller is learned, and corrected for, in the output of the linearizing controller. the result is a controller that \textit{exactly} linearizes the relationship between the inputs and outputs of the system. This technique is implemented and used here for the purpose of learning the dynamics of a 4-wheeled motor vehicle. The resulting linearized dynamics are used to design a path planning and control algorithm based on a Linear Quadratic Regulator, that outputs a sequence of states for the linearized system to follow as well as the sequence of controls needed to achieve those states. While this is sufficient for the exact linearization case, the approximately learned controller will have extra elements to its dynamics that will cause the vehicle to diverge from the optimal path and requiring a more active feedback scheme through the use of the continuous time Linear Quadratic Regulator that calculates appropriate control inputs at every time step, instead of a-priori at the beginning of an episode. We further posit an extension to the method of learning feedback linearization by the introduction a preprocessing network; trained in a supervised manor with data collected from the environment that we use to model the car dynamics in the approximate learned controller case. The purpose of the preprocessing network is to convert the exactly linearized outputs of our learned feedback linearization controller into the inputs required by the environment. Specifically we seek to convert the acceleration value output by our linearizing controller into steering and brake signals that are put into the 'CarRacing-v0'. As of the writing of this report, the results of this work are as follows: a custom gym environment has been created to test the exactly linearizing case of model mismatch in two instances of the same set of equations to test and verify the performance of the LFBL method in a clearly defined case. Further, an initial implementation of the preprocessing network has been trained using data from the target racing environment, but we have been unable to begin experieents on the complete control scheme because of difficulties in the aforementioned exactly linearizing case.   

There are multiple next steps to be performed in continuing this work beyond the progress reported here. Namely, once a hyperparameter sweep is performed with a good initial seed, we expect that a much better performing agent will be learned, we will then begin experiments with the true dynamics in the aforementioned car racing environment. beyond that, once we can demonstrate some basic degree of performance in the car racing environment we can begin tuning the path planning algorithm to optimise more closely for speed and general performance. This would likely involve the iterative use of the path planner to form an MPC-like approach.

\section*{Acknowledgment}
The Authors would like to thank Tyler Westenbroek and Valmik Prabhu for their help in planning experiements and working with the learning feedback linearization method.
\printbibliography

\end{document}